\documentclass[12pt]{article}
\makeindex
\usepackage{epsfig,amsmath,a4wide,amssymb,amsfonts,amsbsy,xy,latexsym,epsf,pstricks,verbatim,times}
\xyoption{all}
\usepackage[french]{babel}

\let\set\mathbb
\def\<<{\leavevmode
  \raise0.28ex\hbox{$\scriptscriptstyle\langle\!\langle$}\nobreak
  \hskip -.6pt plus.3pt minus.2pt\,}
\def\>>{\,\nobreak\hskip -.6pt plus.3pt minus.2pt
  \raise0.28ex\hbox{$\scriptscriptstyle\rangle\!\rangle$}}

\def\dx{{d_x}}
\def\dy{{d_y}}

\def\AA{{\set A}}

\def\Hom{\mathop{\rm{Hom}}\nolimits }

\def\Ker{\mathop{\rm{Ker}}\nolimits }
\def\Gal{\mathop{\rm{Gal}}\nolimits }
\def\Pic{\mathop{\rm{Pic}}\nolimits }

\def\End{\mathop{\rm End }}
\def\FF{{\set F}}

\def\Fp{{\FF _p}}

\def\Fq{{\FF _q}}

\def\bK{{\bf K  }}
\def\bL{{\bf L  }}

\def\bG{{\bf  G }}
\def\bH{{\bf H  }}
\def\oG{{\oplus_\bG}}
\def\oGm{{{\oplus_\bG}_m}}
\def\oGa{{{\oplus_\bG}_a}}
\def\oH{{\oplus_\bH}}
\def\zG{{0_\bG}}
\def\zH{{0_\bH}}

\def\PU{{{\set P}^1}}
\def\AA{{\set A}}
\def\AU{{{\set A}^1}}

\def\ZZ{{\set Z}}
\def\agot{{\mathfrak a}}
\def\igot{{\mathfrak i}}

\def\bgot{{\mathfrak b}}
\def\cA{{\cal A}}
\def\cB{{\cal B}}
\def\cC{{\cal C}}

\def\cF{{\cal F}}
\def\cG{{\cal G}}
\def\cH{{\cal H}}
\def\cI{{\cal I}}

\def\cL{{\cal L}}

\def\cS{{\cal S}}

\def\cU{{\cal U}}

\def\cgot{{\mathfrak  c}}

\def\mmu{{\set \mu}}

\begin{document}
\author{J.-M. Couveignes\thanks{Institut de  Math\'ematiques  de Toulouse,
Universit\'e de Toulouse et CNRS}}
\title{Bases invariantes de friabilit{\'e}\thanks{Cette recherche est soutenue par la
D{\'e}l{\'e}gation G{\'e}n{\'e}rale pour l'Armement,
Centre \'Electronique de l'Armement.}}

\maketitle

\bibliographystyle{plain}

\begin{abstract}

Cette note r{\'e}pond {\`a} une question pos{\'e}e
par Reynald Lercier {\`a} propos des algorithmes de calcul de logarithmes discrets.

\'Etant donn{\'e} un corps r{\'e}siduel fini $k$, on recherche  une  base de
friabilit{\'e} de $k^*$
qui soit stabilis{\'e}e par l'action du groupe des automorphismes de $k$.
Nous construisons des repr{\'e}sentations originales
de
certains corps finis, qui admettent de telles bases.
Ce travail vise {\`a} perfectionner  les algorithmes de calcul
du logarithme discret. On traite le cas de la codimension un  (crible lin{\'e}aire) et de la codimension
deux (crible alg{\'e}brique). 
\end{abstract}

\tableofcontents

\section{Pr{\'e}sentation}\label{section:presentation}

Dans ce travail nous {\'e}tudions l'existence de bases de friabilit{\'e} invariantes
par l'action de Galois pour les corps finis. On sait que de  telles bases 
acc{\'e}l{\`e}rent le calcul
des logarithmes discrets. Nous rappelons cette observation  de Joux et Lercier dans la section
\ref{section:introduction} et donnons un premier exemple dans
la section \ref{section:exemple}. 
Dans la section \ref{section:KAS} nous rappelons les rudiments des th{\'e}ories
de Kummer et Artin-Schreier, qui produisent tous les exemples utilis{\'e}s 
 {\`a} ce jour
de bases Galois invariantes de friabilit{\'e}.
Nous montrons dans la section \ref{section:espaces} que seules les extensions
de Kummer et Artin-Schreier admettent des drapeaux  d'espaces lin{\'e}aires
invariants 
par l'action de Galois. Dans la section \ref{section:isog} nous d{\'e}crivons le
cadre, plus g{\'e}n{\'e}ral, de notre {\'e}tude : la sp{\'e}cialisation d'isog{\'e}nies entre
groupes alg{\'e}briques, et nous en d{\'e}duisons un premier exemple
nouveau de base invariante de friabilit{\'e} dans la section \ref{section:tore}. 
Nous montrons dans la section \ref{section:ell} que les isog{\'e}nies entre
courbes elliptiques
produisent une vari{\'e}t{\'e}  consid{\'e}rable de bases invariantes de friabilit{\'e} lorque
le degr{\'e} du corps n'est pas trop grand.

Dans la section \ref{section:JL} nous rappelons le principe des algorithmes
de cribles rapides, tels que le crible alg{\'e}brique. Nous montrons dans la
section
\ref{section:carre} que notre approche est compatible avec ces principes.
Nous concluons par quelques questions
et remarques sur  les
possibilit{\'e}s 
et les limites de notre m{\'e}thode.

\section{Une question soulev{\'e}e  par Joux et  Lercier}\label{section:introduction}

Rappelons le principe d'un algorithme simple pour calculer des 
logarithmes discrets dans le groupe
multiplicatif
d'un corps fini $\Fq$ avec $q=p^d$. 

Le corps fini $\Fq$  est vu comme corps r{\'e}siduel
$k=\Fp[X]/A(X)$ avec  $A(X)\in \Fp[X]$ polyn{\^o}me unitaire irr{\'e}ductible
de degr{\'e} $d$. On note
$x=X\bmod A(X)$.

Si $n$ est un entier tel que $0\le n\le d-1$ on note $L_n\subset \Fq$ le
$\Fp$-espace
vectoriel engendr{\'e} par $1$, $x$, \dots, $x^{n}$. 

Ainsi $L_0=\Fp\subset L_1\subset \ldots \subset L_{d-1}=\Fq$ et 
$L_a\times L_b\subset L_{a+b}$ si $a+b\le n-1$.

On  construit, par divers moyens, 
des relations multiplicatives entre  {\'e}l{\'e}ments
de $L_\kappa$ o{\`u} $\kappa$ est 
entier
$\kappa$ bien choisi. Par exemple, pour $\kappa=1$, les relations recherch{\'e}es sont   de la forme

\begin{equation}\label{eqn:prod}
\prod_{1\le i \le I} (a_i+b_ix)^{e_i}=1\in \Fq
\end{equation}
\noindent o{\`u} les $a_i$ et $b_i$ sont dans $\Fp$.

On accumule  de telles relations jusqu'{\`a} obtenir une base du
$\ZZ$-module des relations entre les {\'e}l{\'e}ments de $L_\kappa$.

Comment trouve-t-on des relations de type \ref{eqn:prod} ?
Supposons encore que $\kappa=1$.
La forme la plus simple du crible choisit 
des triplets $(a_i,b_i,e_i)$ au hasard et calcule
le reste $r(X)$ de la division euclidienne de $\prod_i(a_i+b_iX)^{e_i}$
par $A(X)$. Donc 

$$r(X)\equiv \prod_i(a_i+b_iX)^{e_i}\bmod A(X)$$
\noindent o{\`u} $r(X)$ est un polyn{\^o}me plus ou moins al{\'e}atoire
 de degr{\'e} $\le d-1$.

On esp{\`e}re que $r(X)$ se d{\'e}compose en produit de polyn{\^o}mes
de degr{\'e} plus
petit que $\kappa=1$. Donc $r(X)=\prod_j(u_j+v_jX)^{f_j}$ et on obtient
la relation

$$\prod_{ i } (a_i+b_ix)^{e_i}\prod_j(u_j+v_jx)^{-f_j}=1$$
\noindent qui est bien du type cherch{\'e}.

On dit que $L_\kappa$ est  la base de friabilit{\'e}.

Joux et Lercier  notent  dans \cite{JouxLercier2} que s'il existe un automorphisme
$\agot$ de $\Fq$ tel que $\agot (x)= ux+v$ avec $u$, $v \in \Fp$,
alors l'action de $\agot$ sur l'{\'e}quation \ref{eqn:prod} produit une
autre {\'e}quation du m{\^e}me type.

Comme l'efficacit{\'e} des algorithmes de calcul du logarithme 
discret  d{\'e}pend du nombre
d'{\'e}qua\-tions de type \ref{eqn:prod}  que l'on peut produire en un temps
donn{\'e}, on souhaite savoir quand de tels automorphismes providentiels existent.

On se demande aussi comment g{\'e}n{\'e}raliser cette observation.

Notons que $\agot$ n'agit pas seulement sur les {\'e}quations (produits) mais aussi sur les
``inconnues'' ou pour mieux dire sur les facteurs  $a_i+b_ix$. Aussi, plut{\^o}t
que d'augmenter le nombre d'{\'e}quations, on peut dire que l'action de $\agot$
permet
de diminuer le nombre d'inconnues (ou de facteurs dans la base de friabilit{\'e}).

En effet si $\agot$ est la puissance $\alpha$-i{\`e}me du Frobenius on obtient la
relation {\it gratuite}

\begin{equation}\label{eqn:free}
\agot(x)=x^{p^\alpha}=ux+v. 
\end{equation}

On peut donc retirer $ux+v$ de la base de friabilit{\'e} et le remplacer partout
par $x^{p^\alpha}$.

Ainsi, on ne conserve qu'un repr{\'e}sentant par orbite de l'action de Galois sur
$L_\kappa$.
Et la taille du syst{\`e}me lin{\'e}aire {\`a} r{\'e}soudre s'en trouve
divis{\'e}e par l'ordre du groupe engendr{\'e} par $\agot$. Si $\agot$
engendre le groupe de Galois de $\Fq/\Fp$ alors on  a divis{\'e}
le nombre d'inconnues  par
$d$, 
le degr{\'e} du corps fini $\Fq$.

Notre pr{\'e}occupation dans ce texte est  de chercher des mod{\`e}les
pour les corps finis, dans lesquels les automorphismes respectent la forme
particuli{\`e}re de certains {\'e}l{\'e}ments o{\`u} de certaines formules.

Par exemple, si le corps fini est pr{\'e}sent{\'e} comme ci-dessus, les
{\'e}l{\'e}ments sont donn{\'e}s comme des polyn{\^o}mes en le g{\'e}n{\'e}rateur
 $x$.   Tout {\'e}l{\'e}ment $z$ du corps fini a un degr{\'e} : c'est le plus
 petit entier $k$ tel que $z\in L_k$. Le degr{\'e}
de  $a_0+a_1x+\dots +a_kx^k$ est donc $k$ pourvu que $0\le k<d$ et $a_k\not =
0$.

Le degr{\'e} est sous-additif $\deg(z\times t)\le \deg(z)+\deg(t)$.

La question pos{\'e}e revient {\`a} se  demander si ce ``degr{\'e}'' 
est pr{\'e}serv{\'e} par les automorphismes de $\Fq$.

On observera que l'int{\'e}r{\^e}t de la fonction degr{\'e} sur $\Fq$ dans le cadre
des algorithmes de cribles tient aux propri{\'e}t{\'e}s suivantes :

\begin{itemize}
\item  le degr{\'e} est sous-additif (et m{\^e}me assez souvent il
est additif) : le degr{\'e} du produit de deux {\'e}l{\'e}ments est la somme
des degr{\'e}s des deux facteurs, pourvu que cette somme soit $<d$.
\item le degr{\'e} permet de repartir agr{\'e}ablement les {\'e}l{\'e}ments de $\Fq$ : 
il y a  $q^n$ {\'e}l{\'e}ments de degr{\'e} $<n$ si $n\le d$.
\item on dispose d'un algorithme de factorisation qui permet de d{\'e}composer ais{\'e}ment
  certains {\'e}l{\'e}ments de $L_{d-1}=\Fq$ en produits d'{\'e}l{\'e}ments de degr{\'e} plus
  petit qu'un  $\kappa$ donn{\'e}. La densit{\'e} dans  $\Fq$  de ces {\'e}l{\'e}ments (dits $\kappa$-friables) n'est pas trop faible.
\end{itemize}

Dans cet article  on part {\`a} la recherche de  fonctions 
``degr{\'e}'' sur les corps finis avec une exigence suppl{\'e}mentaire :
on veut que le degr{\'e} soit invariant par action de Galois.

\section{Un premier exemple}\label{section:exemple}

Voici un premier exemple donn{\'e} par Joux et  Lercier :

Ici $p=43$ et $d=6$ donc $q=43^6$ et on choisit $A(X)=X^6-3$ qui est
bien irr{\'e}ductible dans $\FF_{43}[X]$. Donc $\Fq$ est repr{\'e}sent{\'e}
comme corps r{\'e}siduel 
$k=\FF_{43}[X]/X^6-3$. 

On v{\'e}rifie que 
$p=43$ est congru {\`a} $1$ modulo
$d=6$ donc 

$$\phi(x)= x^{43} = (x^6)^7\times x=3^7x=\zeta_6 x$$
\noindent o{\`u} $\zeta_6=3^7=37\bmod 43$ est une racine sixi{\`e}me primitive
de l'unit{\'e}.

Le Frobenius $\phi$ engendre bien s{\^u}r tout le groupe de Galois. On peut
donc diviser par $6$ la taille de la base de friabilit{\'e}.

Dans le deuxi{\`e}me exemple fourni par Joux et Lercier (issu de XTR de type T30) on
a $p=370801$ et $d=30$ avec $A(X)=X^{30}-17$.
Cette fois $p$ est congru {\`a} $1$ modulo $d=30$  donc

$$\phi(x)=x^p=x^{30\times 12360}\times x=\zeta_{30}x$$
\noindent avec $\zeta_{30}=17^{12360}\bmod p=172960\bmod p$.

Cette fois, on peut diviser par $30$ le nombre d'inconnues.

On est ici dans le cadre de la th{\'e}orie de Kummer. Nous donnerons
donc quelques rappels sur cette th{\'e}orie, qui classifie les extensions
cycliques de $\Fp$ de 
degr{\'e} $d$ divisant $p-1$. La th{\'e}orie d'Artin-Schreier
est le pendant de la th{\'e}orie de Kummer pour les $p$-extensions cycliques en 
caract{\'e}ristique $p$ et nous la pr{\'e}senterons aussi.
Nous allons buter tr{\`e}s vite sur les limitations de
ces deux th{\'e}ories.

Il sera temps alors de consid{\'e}rer la situation plus g{\'e}n{\'e}rale d'un groupe
alg{\'e}brique muni d'un automorphisme rationnel d'ordre fini.

\section{Th{\'e}ories de Kummer et Artin-Schreier}\label{section:KAS}

Il s'agit de classifier les extensions cycliques de degr{\'e} $d$
d'un corps $\bK$ 
de caract{\'e}ristique $p$ dans les
deux cas les   plus simples : 

\begin{itemize}
\item Kummer : 
si $p$ est premier {\`a} $d$
et $\bK$ contient une racine primitive $d$-i{\`e}me de l'unit{\'e};
\item Artin-Schreier : si $d=p$.
\end{itemize}
\vskip 0.3cm

Selon la th{\'e}orie de Kummer, si $p$ est premier {\`a} $d$ 
et $\bK$ contient une racine primitive
de l'unit{\'e}, alors les extensions cycliques de degr{\'e} $d$ 
sont radicielles. Elles se construisents avec des racines.

On consid{\`e}re $r$ un {\'e}l{\'e}ment du  groupe $\bK^*/(\bK^*)^d$ (que l'on identifie
avec un repr{\'e}sentant dans $\bK^*$) et on lui associe
le corps $\bL=\bK(r^\frac{1}{d})$. 

Cette expression  sous-entend que $\bK$ est plong{\'e} dans une cl{\^o}ture 
alg{\'e}brique
$\bar \bK$ et $r^\frac{1}{d}$ est l'une quelconque des racines de l'{\'e}quation
$X^d=r$ dans $\bar\bK$.

On observe 
que l'application $x\mapsto x^d$ d{\'e}finit un {\'e}pimorphisme de
groupe de $\bar\bK^*$ multiplicatif sur lui-m{\^e}me. Le noyau de cet {\'e}pimorphisme
est le groupe des racines $d$-i{\`e}mes de l'unit{\'e}.
Les racines $r^\frac{1}{d}$ ne sont  que les ant{\'e}c{\'e}dents de
$r$ par cet {\'e}pimorphisme.

Le corps $\bK(r^\frac{1}{d})$ n'est pas toujours isomorphe {\`a} l'alg{\`e}bre
$\bK[X]/X^d-r$. Il l'est lorsque $r$ est d'ordre $d$ dans le groupe
$\bK^*/(\bK^*)^d$.

\`A l'extr{\`e}me oppos{\'e}, si $r$ est dans $(\bK^*)^d$
alors $\bK[X]/X^d-r$   est que le produit de $d$ corps isomorphes {\`a} $\bK$.

Revenons au cas o{\`u} $r$ est d'ordre $d$. L'extension $\bL/\bK$ de degr{\'e} $d$ 
est galoisienne car posant $s=r^\frac{1}{d}$ il vient 

$$X^d-r=(X-s)(X-s\zeta_{d})(X-s\zeta_{d}^2)\dots(X-s\zeta_{d}^{d-1})$$
\noindent o{\`u} $\zeta_d$ est une racine primitive $d$-i{\`e}me de
l'unit{\'e}.

Le groupe de Galois de $\bL/\bK$ 
est form{\'e} des transformations de la forme

$$\agot_k : s\mapsto s\zeta_d^k$$
\noindent et l'application
$k\mapsto \agot_k$ est un isomorphisme du groupe $\ZZ/d\ZZ$ vers
$\Gal(\bL/\bK)$.

Si l'on veut {\'e}viter de distinguer une infinit{\'e} de cas, selon que $r$ est
d'ordre petit ou grand dans  $\bK^*/(\bK^*)^d$, on proc{\`e}de comme
dans Bourbaki \cite[A V.84]{bourbaki}.

Plut{\^o}t que de prendre un {\'e}l{\'e}ment de $\bK^*/(\bK^*)^d$ on choisit 
un sous-groupe $H$ de $\bK^*$ contenant $(\bK^*)^d$ et on forme
l'extension $\bK(H^\frac{1}{d})$ en prenant toutes les racines $d$-i{\`e}mes
des {\'e}l{\'e}ments de $H$.

\`A tout {\'e}l{\'e}ment $\agot$ de $\Gal(\bK(H^\frac{1}{d})/\bK)$ 
on associe alors
un homomorphisme $\kappa(\agot)$ de $H/(\bK^*)^d$ vers le groupe $\mmu_d$
des racines  $d$-i{\`e}mes de l'unit{\'e}. L'homomorphisme 
$\kappa (\agot)$  est d{\'e}fini par 

$$\kappa(\agot) : \theta \mapsto
\frac{\agot(\theta^\frac{1}{d})}{\theta^\frac{1}{d}}$$
\noindent o{\`u} $\theta^\frac{1}{d}$ est l'une des racines
$d$-i{\`e}mes de $\theta$ (mais on doit bien s{\^u}r prendre la m{\^e}me 
au num{\'e}rateur et au d{\'e}nominateur !)

La correspondance $\agot \mapsto \kappa (\agot)$ est
un isomorphisme du groupe de Galois 

$$\Gal(\bK(H^\frac{1}{d})/\bK)$$
\noindent 
vers le groupe des homomorphismes $\Hom(H/(\bK^*)^d, \mmu_d)$.

Cela revient {\`a} caract{\'e}riser un automorphisme $\agot$
par la mani{\`e}re dont  il agit sur certains radicaux.

Cette pr{\'e}sentation de la th{\'e}orie de Kummer construit les extensions
ab{\'e}liennes de $\bK$ d'exposant divisant $d$.

Dans le cas qui nous int{\'e}resse  le corps $\bK=\Fq$ est fini. 
Tout sous groupe $H$ de $\bK^*$ est cyclique. Pour avoir
$\mmu_d$ dans $\bK$ on doit supposer que $d$ divise $q-1$.
On note $q-1=md$. Le groupe $(\bK^*)^d$ a pour
cardinal $m$. Le quotient $\bK^*/(\bK^*)^d$
est cyclique d'ordre  $d$ donc il est naturel de choisir $H=\bK^*$
(on ne peut mieux faire.)

On en d{\'e}duit qu'il existe une unique extension cyclique $\bL$
de degr{\'e}  $d$,
 engendr{\'e}e
par une racine $d$-i{\`e}me d'un g{\'e}n{\'e}rateur $r$ de $\bK^*$.

Soit donc $s=r^\frac{1}{d}$ et $\bL=\bK(s)$. 
Le  groupe de Galois $\Gal(\bL/\bK)$ est engendr{\'e} par le Frobenius
$\phi$ et l'action de $\phi$ sur $s$ est donn{\'e}e par
$\phi(s)=s^{q}$ donc

$$\frac{\phi(s)}{s}=s^{q-1}=\zeta=r^m$$
\noindent o{\`u} $\zeta$ est une racine $d$-i{\`e}me de l'unit{\'e} qui d{\'e}pend 
de $r$. La correspondance $r\mapsto \zeta$ est un isomorphisme
du groupe $\bK^*/(\bK^*)^d$  sur le groupe $\mmu_d$ qui n'est autre
que l'exponentiation par $m$.

Revenant au premier exemple on a $q=p=43$, $p-1=42$, $d=6$, $m=7$,
$r=3$ et $\frac{\phi(s)}{s}=r^m=3^7\bmod 43$.

On voit imm{\'e}diatement les limites de cette construction : elle requiert
la pr{\'e}sence de racines $d$-i{\`e}mes primitives de  l'unit{\'e} dans $\bK$.

Si ces racines font d{\'e}faut, on doit recourir {\`a} d'autres moyens
pour construire des extension de corps cycliques. Les 
automorphismes des extensions construites par ces m{\'e}thodes
plus g{\'e}n{\'e}rales ne semblent pas souffrir une pr{\'e}sentation aussi
simple que dans la th{\'e}orie de Kummer. 

On peut par exemple  passer par une extension auxiliaire
$\bK'=\bK(\zeta_d)$  de $\bK$,  qui peut {\^e}tre h{\'e}las tr{\`e}s grande. On applique
alors la th{\'e}orie de Kummer {\`a} cette grosse extension et
on obtient une extension $\bL'/\bK'$ cyclique de degr{\'e} $d$.
  La descente de cette extension se fait par des moyens alg{\'e}briques (resolvantes)
peu compatibles avec les  exigences formul{\'e}es dans la section \ref{section:introduction}.
On pourra voir \cite[Chapitre III.4]{matzat}. Nous  n'explorerons donc
pas cette piste.

On r{\'e}sume maintenant la th{\'e}orie d'Artin-Schreier.

Selon cette  th{\'e}orie, si $p$  est la caract{\'e}ristique de $\bK$
 alors toute  extension cyclique de degr{\'e} $p$ 
est engendr{\'e}e par les racines d'un polyn{\^o}me de la forme

$$X^p-X-a=\wp(X)-a=0$$
\noindent o{\`u} $a\in \bK$ et o{\`u} l'expression $\wp(X)=X^p-X$ semble jouer ici un r{\^o}le
assez comparable {\`a} celui de $X^n$ dans la th{\'e}orie de Kummer.

On observe en effet
que l'application $x\mapsto \wp(x)$ d{\'e}finit un {\'e}pimorphisme de
groupe de $\bar\bK$ additif sur lui-m{\^e}me. Le noyau de cet {\'e}pimorphisme
est le groupe additif du corps premier $\Fp\subset \bar\bK$.

On consid{\`e}re $a$ un {\'e}l{\'e}ment du  groupe additif $\bK/\wp(\bK)$ (que l'on identifie
avec un repr{\'e}sentant dans $\bK$) et on lui associe
le corps $\bL=\bK(\wp^{-1}(a))$. 

Ici encore  on sous-entend que $\bK$ est plong{\'e} dans une cl{\^o}ture 
alg{\'e}brique
$\bar \bK$.  Alors $\bL$ est le sous-corps de $\bar \bK$ engendr{\'e} par
$\bK$ et l'une quelconque des racines de l'{\'e}quation $\wp(X)=a$. Comme
deux racines diff{\'e}rent d'un {\'e}l{\'e}ment du corps primitif $\Fp$, il importe
peu de savoir laquelle on a choisie.

Ici encore, le corps $\bK(\wp^{-1}(a))$ n'est pas toujours isomorphe {\`a} l'alg{\`e}bre
$\bK[X]/X^p-X-a$. Il l'est lorsque $a$ est non nul dans $\bK/\wp(\bK)$.
Sinon,  $\bK[X]/X^p-X-a$   est le produit de $p$ corps isomorphes {\`a} $\bK$.

On suppose donc que  $a$ est non-nul donc d'ordre $p$ dans $\bK/\wp(\bK)$.
L'extension $\bL/\bK$ de degr{\'e} $p$ 
est galoisienne car posant $b=\wp^{-1}(a)$ il vient 

$$X^p-X-a=(X-b)(X-b-1)(X-b-2)\dots(X-b-(p-1)).$$

Le groupe de Galois est form{\'e} des transformations de la forme

$$\agot_k : b\mapsto b+k$$
\noindent et l'application
$k\mapsto \agot_k$ est un isomorphisme du groupe $\ZZ/p\ZZ$ vers
$\Gal(\bL/\bK)$.

Si l'on veut maintenant construire toutes les
extensions ab{\'e}liennes de $\bK$ d'exposant $p$,
on suit  Bourbaki \cite[A V.88]{bourbaki}. On consid{\`e}re   un 
sous-groupe $H$ de $(\bK,+)$ contenant $\wp(\bK)$ et on forme
l'extension $\bK(\wp^{-1}(H))$.

\`A tout {\'e}l{\'e}ment $\agot$ de $\Gal(\bK(\wp^{-1}(H))/\bK)$ 
on associe alors
un homomorphisme $\kappa(\agot)$ de $H/\wp(\bK)$ vers le groupe 
additif $\Fp$
du corps premier. L'homomorphisme 
$\kappa (\agot)$  est d{\'e}fini par 

$$\kappa(\agot) : \theta \mapsto
\agot(\wp^{-1}(\theta))- \wp^{-1}(\theta)$$
\noindent o{\`u} $\wp^{-1}(\theta)$ est l'un des ant{\'e}c{\'e}dents
de $\theta$ par $\wp$ (et bien s{\^u}r on doit prendre le m{\^e}me dans
le premier et dans le second terme de la diff{\'e}rence.)

La correspondance $\agot \mapsto \kappa (\agot)$ est
un isomorphisme du groupe de Galois 

$$\Gal(\bK(\wp^{-1}(H))/\bK)$$
\noindent 
vers le groupe des homomorphismes $\Hom(H/\wp(\bK), \Fp)$.

Dans le cas qui nous int{\'e}resse  le corps $\bK=\Fq$ est fini de
caract{\'e}ristique 
$p$.  On pose $q=p^f$.

Le morphisme 
$$\wp : \Fq \rightarrow \Fq$$
\noindent  a pour noyau
$\Fp$ donc le quotient $\Fq/\wp(\Fq)$ est d'ordre $p$.

Il existe donc une seule extension $\bL$ de degr{\'e} $p$ de $\Fq$ et elle
est engendr{\'e}e par $b=\wp^{-1}(a)$ avec $a\in \Fq-\wp(\Fq)$.

Le  groupe de Galois $\Gal(\bL/\bK)$ est engendr{\'e} 
par le Frobenius $\phi$ et 
$\phi(b)-b$ appartient {\`a} $\Fp$.  La correspondance $a\mapsto \phi(b)-b$ 
est un isomorphisme
du groupe $\bK/\wp(\bK)$  sur le groupe additif  $\Fp$.

On se demande s'il est possible de rendre plus explicite 
cet isomorphisme.

On a $\phi(b)=b^q$ o{\`u} $q=p^f$ est le cardinal de $\bK=\Fq$. Donc 

$$\phi(b)-b=b^q-b=(b^p)^{p^{f-1}}-b=(b+a)^{p^{f-1}}-b$$
\noindent 
car $\wp(b)=b^p-b=a$.

Donc $b^{p^{f}}-b=b^{p^{f-1}}-b +a^{p^{f-1}}$. En it{\'e}rant il
vient 

$$\phi(b)-b=b^{p^{f}}-b=a+a^p+a^{p^2}+\cdots+a^{p^{f-1}}.$$

Ainsi l'isomorphisme du groupe $\bK/\wp(\bK)$  sur le groupe additif  $\Fp$
n'est autre que la trace absolue.

{\bf Exemple} : on choisit $p=7$ et $f=1$
donc $q=7$. La trace absolue de $1$ est $1$
donc on pose $\bK=\FF_7$ et $A(X)=X^7-X-1$ et on construit
$\bL=\FF_{7^7}=\FF_7[X]/A(X)$.  On note $x=X\bmod A(X)$.

On a $\phi(x)=x+1$.

\section{Sous-espaces lin{\'e}aires invariants d'une
extension cyclique}\label{section:espaces}

On rappelle que la question pos{\'e}e dans la section
\ref{section:introduction} revient {\`a} se demander
s'il existe des automorphismes qui respectent une certaine
base de friabilit{\'e}. 

On a vu dans l'introduction que les bases de friabilit{\'e} sont form{\'e}es
ordinairement {\`a} l'aide d'un drapeau d'espaces vectoriels.

On se demande donc si, pour une extension cyclique
 $\bL/\bK$ donn{\'e}e, il existe des $\bK$-sous-espaces vectoriels
de $\bL$ invariants par le groupe de Galois de $\bL/\bK$.

Supposons que $\bL=\bK[X]/X^d-r$ est une extension de Kummer
et pour tout entier $k$ entre $0$ et $d-1$ notons 

$$L_k=\bK\oplus \bK x\oplus \cdots\oplus \bK x^k$$
\noindent  le $\bK$-sous-espace vectoriel engendr{\'e} par les
$k+1$ premi{\`e}res puissances de $x=X\bmod X^d-r$.

Les $L_k$ sont globalement invariants par l'action de Galois car
si $\agot$ est un 

\noindent
$\bK$-automorphisme de $\bL$  alors
il existe une racine $d$-i{\`e}me de l'unit{\'e} $\zeta \in \bK$
telle que  

$$\agot(x)=\zeta x$$
\noindent  et  $\agot(x^k)=\zeta^k x^k$.

On a donc un drapeau de $\bK$-espaces vectoriels 

$$\bK=L_0\subset L_1\subset \dots \subset L_{d-1}=\bL$$
\noindent qui est respect{\'e} par l'action de Galois. 
Donc le ``degr{\'e}'' est respect{\'e} par cette action.

C'est tr{\`e}s exactement ce qui se produit dans les deux exemples
de la section \ref{section:introduction}
 : le ``degr{\'e}'' des {\'e}l{\'e}ments du corps fini est respect{\'e} 
par l'action de Galois. Donc si la base de friabilit{\'e} est constitu{\'e}e
par tous les polyn{\^o}mes irr{\'e}ductibles
de degr{\'e} $\le \kappa$ alors elle est globalement 
invariante par l'action de Galois.

Supposons maintenant que $\bL=\bK[X]/X^p-X-a$ est une extension 
d'Artin-Schreier
et pour tout entier $k$ entre $0$ et $p-1$ notons 

$$L_k=\bK\oplus \bK x\oplus \cdots\oplus \bK x^k$$
\noindent  le $\bK$-sous-espace vectoriel engendr{\'e} par les
$k+1$ premi{\`e}res puissances de $x=X\bmod X^p-X-a$.

Les $L_k$ sont globalement invariants par l'action de Galois car
si $\agot$ est un 
$\bK$-automorphisme de $\bL$  alors
il existe une constante   $c \in \Fp$
telle que  

$$\agot(x)=x+c$$
\noindent  et  donc

$$\agot(x^k)=(x+c)^k=\sum_{0\le \ell\le k}\left( \begin{array}{c}   
k\\ \ell
\end{array}  \right)c^{k-\ell}x^\ell.$$

On a donc  encore un drapeau de $\bK$-espaces vectoriels 

$$\bK=L_0\subset L_1\subset \dots \subset L_{p-1}=\bL$$
\noindent qui est respect{\'e} par l'action de Galois.

Attention ! Cette fois, l'action de Galois n'est pas diagonale mais
seulement triangulaire.

Notons que pour les extensions de degr{\'e} une puissance de $p$, la
 th{\'e}orie dite de Witt-Artin-Schreier  g{\'e}n{\'e}ralise
la th{\'e}orie d'Artin-Schreier. Et elle produit aussi un drapeau
invariant  de sous-espaces vectoriels. On trouve
au  d{\'e}but de la th{\`e}se de Lara Thomas \cite{lara}  des r{\'e}f{\'e}rences et une
pr{\'e}sentation de cette th{\'e}orie.

On peut se demander si des drapeaux d'espaces lin{\'e}aires
invariants existent dans d'autres
cas. 

On suppose que $\bL/\bK$ est une extension cyclique de degr{\'e} $d$ fini et
premier
{\`a} la caract{\'e}ristique $p$. Soit $\phi$ un g{\'e}n{\'e}rateur
de  $C=<\phi>=\Gal(\bL/\bK)$
le groupe de Galois. D'apr{\`e}s le th{\'e}or{\`e}me de la base normale 
\cite[Theorem 13.1.]{langalgebra} il existe un {\'e}l{\'e}ment $w$  de $\bL$
tel que 

$$(w,\phi(w),\phi^2(w),\dots,\phi^{d-1}(w))$$
\noindent  soit une $\bK$-base de
$\bL$.

On en d{\'e}duit que $\bL$ muni de l'action de $C$ est la repr{\'e}sentation
r{\'e}guli{\`e}re de ce groupe cyclique d'ordre $d$ sur le corps
$\bK$.

Comme l'ordre $d$ du groupe $C$ est premier {\`a} la  caract{\'e}ristique $p$
de $\bK$, l'anneau $\bK[C]$ est semi-simple d'apr{\`e}s le th{\'e}or{\`e}me
de Maschke \cite[Theorem 1.2.]{langalgebra}. Cela
signifie que toute repr{\'e}sentation est somme directe de repr{\'e}sentations
irr{\'e}ductibles. Autrement dit ``tout se passe comme en caract{\'e}ristique
z{\'e}ro''.

Le polyn{\^o}me caract{\'e}ristique de $\phi$ sur le $\bK$ espace vectoriel
$\bL$ est $X^d-1$ qui est un polyn{\^o}me  s{\'e}parable sur $\bK$.

Il existe pour chaque $\bK$-facteur $f(X)\in \bK[X]$ de $X^n-1$ 
un unique sous-espace irr{\'e}ductible $V_f\subset \bL$ invariant
par $\phi$ et tel que
la restriction de $\phi$ {\`a} $V_f$ admette $f$ comme polyn{\^o}me
caract{\'e}ristique.

Tout sous-espace invariant par $\phi$ est somme directe de quelques
$V_f$ d'apr{\`e}s  le lemme de Schur \cite[Proposition 1.1.]{langalgebra}.

Pour qu'il existe un drapeau complet de sous-espaces invariants
par $\phi$

$$\bK=L_0\subset L_1\subset \dots \subset L_{d-1}=\bL$$
\noindent avec $L_k$ de dimension $k$, on doit avoir uniquement des 
facteur irr{\'e}ductibles de degr{\'e} $1$ dans $X^d-1$.

Alors $\bK$ contient les racines primitives $d$-i{\`e}mes de l'unit{\'e}
et on est dans le cadre de la th{\'e}orie de Kummer.

 Au passage, on a d{\'e}montr{\'e} que le drapeau fourni par 
la th{\'e}orie de Kummer est unique {\`a} permutation pr{\`e}s. 
Plus pr{\'e}cis{\'e}ment,
tout drapeau $\phi$-invariant est d{\'e}termin{\'e} par l'ordre choisi sur les racines
de l'unit{\'e} (et donc sur les $\bK[\phi]$-espaces irr{\'e}ductibles
de $\bL$). Il y a $d!$ tels drapeaux. 

Le drapeaux produits par la th{\'e}orie de Kummer ont une propri{\'e}t{\'e} 
suppl{\'e}mentaire : ils sont de la forme

$$V_1\subset V_1\oplus V_\zeta\subset V_1\oplus V_\zeta\oplus V_{\zeta^2}\subset  \dots \subset 
V_1\oplus V_\zeta\oplus V_{\zeta^2}\oplus\cdots\oplus V_{\zeta^{n-2}}
\subset V_1\oplus V_\zeta\oplus V_{\zeta^2}\oplus\cdots\oplus V_{\zeta^{d-2}}\oplus
V_{\zeta^{d-1}}$$
\noindent o{\`u} $\zeta $ est une racine primitive $d$-i{\`e}me de l'unit{\'e}
et $V_\zeta$ est $V_{X-\zeta}$ l'espace irr{\'e}ductible associ{\'e} au
facteur $X-\zeta$ de $X^d-1$ ou, si l'on pr{\'e}f{\`e}re, l'espace propre
associ{\'e} {\`a} la valeur propre $\zeta$ de $\phi$.

Parmi les $d!$ drapeaux $\phi$-invariants disponibles, il y en a 
$\phi (d)$ qui sont fournis par la th{\'e}orie de Kummer. Ils
correspondent aux $\phi(d)$ racines primitives $d$-i{\`e}mes de l'unit{\'e}.
 Ces derniers drapeaux 
jouissent d'une propri{\'e}t{\'e} multiplicative essentielle pour les
applications envisag{\'e}es :  si $k\ge 0$ et $l\ge 0$
et $k+l\le d-1$ alors 

$$L_k\times L_l\subset L_{k+l}.$$

La conclusion de cette section est donc assez n{\'e}gative. Pour aller
plus loin que la th{\'e}orie de Kummer, il faudra se montrer moins exigeant.

\section{Sp{\'e}cialisation d'isog{\'e}nies
entre groupes alg{\'e}briques}\label{section:isog}

La th{\'e}orie de Kummer et la th{\'e}orie d'Artin-Schreier sont deux
cas particuliers d'une situation plus g{\'e}n{\'e}rale que nous allons 
d{\'e}crire
maintenant et qui nous permettra de construire de nouveaux exemples
d'automorphismes agr{\'e}ables pour les corps finis.

Soit $\bK$ un corps et $\bG$  un groupe alg{\'e}brique commutatif.
Soit $T\subset \bG(\bK)$ un groupe fini de points $\bK$-rationnels
de $\bG$  et soit

$$I : \bG\rightarrow \bH$$
\noindent l'isog{\'e}nie quotient de $\bG$ par $T$.

On note $d$ le cardinal de $T$ qui est aussi le degr{\'e} de $I$.

On suppose qu'il existe un point $\bK$-rationnel $a$ sur
$\bH$ tel que $I^{-1}(a)$ 
soit r{\'e}duit  et irr{\'e}ductible sur $\bK$. Donc  tout point $b$
tel que $I(b)=a$ d{\'e}finit une extension $\bL$ de degr{\'e} $d$
de $\bK$. 

On note $\bL=\bK(b)$ et on observe que l'origine g{\'e}om{\'e}trique
de cette extension fournit des $\bK$-automorphismes de $\bL$.

Soit $t$ un {\'e}l{\'e}ment de $T$ et notons $\oG$ l'addition
dans le groupe alg{\'e}brique $\bG$ et $\oH$ l'addition dans $\bH$.

On note $\zG$ l'{\'e}l{\'e}ment neutre de $\bG$ et $\zH$ celui de
$\bH$.

Le point $t\oG b$ v{\'e}rifie

$$I(t\oG b)=I(t)\oH I(b)=\zH\oH a=a.$$

Donc $t\oG b$ est conjugu{\'e} de $b$ par l'action de Galois et on
obtient tous les conjugu{\'e}s de $b$ en prenant tous les $t$
dans  $T$.

On a donc un isomorphisme 
entre $T$ et $\Gal(\bL/\bK)$, qui {\`a} tout $t\in T$ 
associe l'automorphisme  r{\'e}siduel

$$b\in I^{-1}(a)\mapsto b\oG t.$$

Maintenant,  si les formules g{\'e}om{\'e}triques pour la translation 
$P\mapsto P\oG t$ dans $\bG$ sont simples, on a obtenu une description
agr{\'e}able du groupe de Galois de $\bL$ sur $\bK$.

Nous nous int{\'e}ressons ici aux corps finis. Donc
$\bK=\Fq$. Il suffit alors  de trouver un point $b$ dans $\bH(\Fq)$ tel que
$I^{-1}(b)$ soit $\bK$-irr{\'e}ductible. Cela signifie que
les points g{\'e}om{\'e}triques $b$  dans $I^{-1}(a)$ sont d{\'e}finis sur
$\bL=\FF_{q^d}$ et sur aucune sous-extension.

On illustre ces g{\'e}n{\'e}ralit{\'e}s en revenant aux th{\'e}ories 
de Kummer et Artin-Schreier que l'on revoit ici dans le cadre
plus g{\'e}om{\'e}trique que nous venons d'esquisser.

Pour la th{\'e}orie de Kummer, le groupe alg{\'e}brique sous-jacent est
le groupe multiplicatif $\bG_m$. L'isog{\'e}nie $I$ est la multiplication
par $d$ :

$$I=[d] : \bG_m\rightarrow \bG_m.$$

Le groupe $\bG_m$ est vu comme sous-vari{\'e}t{\'e} de la droite
affine  $\bG_m\subset \AA^1$. Un point $P$ de $\bG_m$ est d{\'e}fini
par une seule coordonn{\'e}e $z$. En fait $\bG_m$ est d{\'e}fini par
l'in{\'e}galit{\'e} $z\not =0$.

L'origine $\zG$ a pour coordonn{\'e}e $z(\zG)=1$. La loi de groupe
alg{\'e}brique est donn{\'e}e par

$$z(P_1\oGm P_2)=z(P_1)\times z(P_2).$$

On a ici $\bH=\bG=\bG_m$ et l'isog{\'e}nie $I$ est d{\'e}crite en termes
des coordonn{\'e}es $z$ par 

$$z(I(P))=z(P)^d.$$

Les points du noyau de $I$ ont pour $z$-coordonn{\'e}es les racines
$d$-i{\`e}mes de l'unit{\'e}. 

L'image r{\'e}ciproque par $I$ d'un point $P$ de $\bG$ est form{\'e}e
de $d$ points g{\'e}om{\'e}triques dont les $z$-coordonn{\'e}es  sont 
les $d$ racines $d$-i{\`e}mes de $z(P)$.

La translation par un {\'e}l{\'e}ment $t$ du noyau de $I$ 

$$P\mapsto P\oGm t$$
\noindent 
s'exprime en terme
de $z$-coordonn{\'e}es  par

$$z(P\oGm t)=z(P)\times \zeta$$
\noindent o{\`u} $\zeta=z(t)$ est la racine $d$-i{\`e}me de l'unit{\'e} associ{\'e}e
par $z$ au point de $d$-torsion $t$.

Pour la th{\'e}orie d'Artin-Schreier, le groupe alg{\'e}brique sous-jacent est
le groupe additif $\bG_a$ sur un corps de caract{\'e}ristique $p$. 

Le groupe $\bG_a$ est identifi{\'e} {\`a} la  droite
affine  $\AA^1$. Un point $P$ de $\bG_a$ est d{\'e}fini
par une seule coordonn{\'e}e $z$. 

L'origine $\zG$ a pour coordonn{\'e}e $z(\zG)=0$. La loi de groupe
alg{\'e}brique est donn{\'e}e par

$$z(P_1\oGa P_2)=z(P_1)+ z(P_2).$$

L'isog{\'e}nie $I$ est l'aplication
s{\'e}parable de degr{\'e} $p$ :

$$\wp : \bG_a\rightarrow \bG_a$$
\noindent d{\'e}crite en termes
des coordonn{\'e}es $z$ par 

$$z(\wp(P))=z(P)^p-z(P).$$

On a encore ici $\bH=\bG$.

Les points du noyau de $\wp$ ont pour $z$-coordonn{\'e}es les 
{\'e}l{\'e}ments du corps premier $\Fp$.

L'image r{\'e}ciproque par $I$ d'un point $P$ de $\bG$ est form{\'e}
de $p$ points g{\'e}om{\'e}triques dont les $z$-coordonn{\'e}es  sont 
les $p$ racines  de l'{\'e}quation $X^p-X=z(P)$.

La translation par un {\'e}l{\'e}ment $t$ du noyau de $I$ 

$$P\mapsto P\oGa t$$
\noindent 
s'exprime en terme
de $z$-coordonn{\'e}es  par

$$z(P\oGa t)=z(P) + c$$
\noindent o{\`u} $c=z(t)\in \Fp$.

\section{Un exemple diff{\'e}rent}\label{section:tore}

On veut appliquer  les g{\'e}n{\'e}ralit{\'e}s de la section pr{\'e}c{\'e}dente
{\`a} divers groupes alg{\'e}briques commutatifs. On devine   que chaque
groupe alg{\'e}brique apportera sa petite contribution {\`a} notre
probl{\`e}me. Cependant, comme on cherche des formules simples 
pour la translation, on imagine que ce sont les groupes
alg{\'e}briques les plus ordinaires qui seront les plus utiles.

On commence donc par les plus familiers des groupes alg{\'e}briques
apr{\`e}s les groupes $\bG_m$ et $\bG_a$ : il s'agit
des tores de dimension $1$.

Soit $\bK$ un corps de caract{\'e}ristique diff{\'e}rente de $2$ et 
$D$ un {\'e}l{\'e}ment non-nul de $\bK$.

Soit $\PU$ la droite projective et $[U,V]$ des
coordonn{\'e}es projectives sur $\PU$. On note

$u=\frac{U}{V}$
\noindent la coordonn{\'e}e affine associ{\'e}e.

Soit $\bG$ l'ouvert de $\PU$ d'in{\'e}quation

$$U^2-DV^2\not =0.$$

On associe {\`a} chaque point $P$ de $\bG$ sa $u$-coordonn{\'e}e
{\'e}ventuellement infinie mais distincte de $\sqrt D$ et $-\sqrt D$.

L'{\'e}l{\'e}ment neutre de $\bG$ est le point $\zG$ de coordonn{\'e}es
$[1,0]$ et de $u$-coordonn{\'e}e $\infty$.

La loi d'addition est d{\'e}finie par

$$u(P_1\oG P_2)=\frac{u(P_1)u(P_2)+D}{u(P_1)+u(P_2)}$$
\noindent et

$$u(\ominus_\bG P_1)=-u(P_1).$$

On suppose d{\'e}sormais que $\bK=\Fq$ est un corps fini et que 
$D\in \Fq^*$ n'est pas un carr{\'e} dans $\Fq$.

Le groupe $\bG(\Fq)$ des points $\Fq$-rationnels est de cardinal
$q+1$ et les valeurs correspondantes de $u$ sont dans 
$\Fq\cup \{\infty \}$.

L'endomorphisme de Frobenius

\begin{eqnarray*}\xymatrix{
\phi : & \bG \ar@{->}[r]   &    \bG\\
&[U,V]\ar@{->}[r]  &[U^q,V^q]
}
\end{eqnarray*}
\noindent se confond avec l'isog{\'e}nie multiplication par
$-q$. En effet, soit $P$ le point de coordonn{\'e}es
projectives $[U,V]$. Les coordonn{\'e}es projectives  de
$R=[q]P$ sont les coordonn{\'e}es dans la base $(1, \sqrt{D})$ de 

$$(U+ V\sqrt{D})^q=U^q-\sqrt{D}V^q$$
\noindent car $D$ n'est pas un carr{\'e} dans $\Fq$.

Donc $R$ a pour coordonn{\'e}es $[U^q,-V^q]$ et c'est 
bien l'inverse de $\phi(P)$.

On se donne alors un entier $d\ge 2$  et on demande
que la $d$-torsion $\bG[d]$ soit $\Fq$-rationnelle.
Il faut  que $d$ divise $q+1$. On pose
 $q+1=md$.

On consid{\`e}re l'isog{\'e}nie $I$ multiplication par
$d$ :

$$I=[d]  : \bG\rightarrow \bG$$
\noindent dont le noyau $\bG[d]$ est cyclique
d'ordre $d$ et d{\'e}compos{\'e} sur $\bK=\Fq$. 

Le quotient $\bG(\Fq)/I(\bG(\Fq))=\bG(\Fq)/\bG(\Fq)^d$
est cyclique de cardinal $d$.

Soit alors $r$ un g{\'e}n{\'e}rateur de $\bG(\Fq)$ et soit
$s$ un ant{\'e}c{\'e}dent de $r$ par $I$. On note
$u(s)$ la $u$-coordonn{\'e}e de $s$ et on pose
$\bL=\bK(u(s))$. C'est une extension de degr{\'e} $d$ de $\bK$.

Le groupe de Galois de $\bL/\bK$  est isomorphe {\`a} 
$\bG[d]$ :  pour tout $\agot \in \Gal (\bL/\bK)$, la diff{\'e}rence
$\agot(s)\ominus_\bG s$ est dans $\bG[d]$ et l'accouplement 

$$(\agot , r)\mapsto  \agot(s)\ominus_\bG s$$
\noindent d{\'e}finit un isomorphisme de 
$\Gal(\bL/\bK)$ vers $\Hom(\bG(\bK)/(\bG(\bK))^d,\bG[n])$.

Ici $\Gal(\bL/\bK)$ est cyclique d'ordre $d$ et engendr{\'e} 
par le Frobenius $\phi$. L'accouplement
$(\phi,r)$ vaut $\phi(s)\ominus_\bG s$.

Ici attention, on se souvient que $\phi(s)=[{-q}]s$  dans $\bG$ donc

\begin{equation}\label{eqn:frob}
(\phi,r)=[-q-1]s=[-m]r.
\end{equation}
\noindent

On a donc une description exacte de l'action de Galois sur
$I^{-1}(r)$.
Elle est donn{\'e}e par une translation du type
$P\mapsto P\oG t$ avec  $t\in \bG[d]$. Si la
coordonn{\'e}e affine de  $t$ est
 $\tau$ et si celle de $P$ est $u$ alors
l'action de la translation sur la coordonn{\'e}e $u$ est donn{\'e}e par

$$u\mapsto \frac{\tau u+D}{u+\tau}$$
\noindent qui est tr{\`e}s agr{\'e}able car c'est une homographie.

On forme le polyn{\^o}me

$$A(X)=\prod_{s\in I^{-1}(r)} (X-u(s))$$
\noindent  annulateur des $u$-coordonn{\'e}es
des ant{\'e}c{\'e}dents de $r$ par $I$.

C'est un polyn{\^o}me de degr{\'e} $d$ {\`a} coefficients dans $\bK=\Fq$.
Il est irr{\'e}ductible dans $\Fq[X]$ car $r$ est un g{\'e}n{\'e}rateur
de $\bG (\Fq)$. Donc on construit $\bL$ comme
$\bK[X]/A(X)$.

Les formules d'exponentiation dans $\bG$ permettent de
donner explicitement le polyn{\^o}me $A(X)$.

On a 

$$(U+\sqrt DV)^d=\sum_{0\le 2k\le d} 
\left(\begin{array}{c}
d\\2k \end{array}\right)
U^{d-2k}V^{2k}D^k+\sqrt D\sum_{1\le 2k+1\le
  d} 
\left(\begin{array}{c}
d\\2k+1 \end{array}\right)
U^{d-2k-1}V^{2k+1}D^k.$$

 Donc

$$u([k]P)=\frac{\sum_{0\le 2k\le d} u(P)^{d-2k}
\left(\begin{array}{c}
d\\2k \end{array}\right)
D^k}{\sum_{1\le 2k+1\le
  d} u(P)^{d-2k-1}
\left(\begin{array}{c}
d\\2k+1 \end{array}\right)
D^k}.$$

Ainsi

$$A(X)=\sum_{0\le 2k\le d} X^{d-2k}
\left(\begin{array}{c}
d\\2k \end{array}\right)
D^k-u(r)\sum_{1\le 2k+1\le
  d} X^{d-2k-1}
\left(\begin{array}{c}
d\\2k+1 \end{array}\right)
D^k.$$

On pose $x=X\bmod A(X)$. Puisque tout {\'e}l{\'e}ment du  groupe de Galois
transforme $x$ en une fraction rationnelle de degr{\'e} $1$
en $x$ il est naturel de d{\'e}finir pour tout entier $k$
tel que $k\ge 0$ et  $k< d$ le sous-ensemble

$$P_k=\{ \frac{a_0+a_1x+a_2x^2+\cdots+a_kx^k}{b_0+b_1x+b_2x^2+\cdots+b_kx^k}
  | (a_0,a_1,\ldots,a_k,b_0,b_1,\ldots,b_k)\in \bK^{2k+2}  \}.$$

On a 

$$\bK=P_0\subset P_1\subset \dots\subset P_{d-1}=\bL$$
\noindent et les $P_k$ sont invariants
par action de Galois.

En outre il est clair que

$$P_k\times P_l\subset P_{k+l}$$
\noindent si $k+l\le d-1$.

Donc on a encore un drapeau de sous-ensembles stables par
l'action de Galois mais ces ensembles ne sont
pas lin{\'e}aires.

Si on d{\'e}finit le ``degr{\'e}'' d'un {\'e}l{\'e}ment de $\bL$ comme
le plus petit $k$ tel que $P_k$ contient cet {\'e}l{\'e}ment,
alors le degr{\'e} est une fonction invariante par l'action
de Galois et sous-additive : 

$$\deg(ab)\le \deg(a)+\deg(b).$$

On voit en outre que le degr{\'e} est compris
entre $0$ et $\lceil \frac{d-2}{2} \rceil $. C'est donc 
une fonction un peu moins fine que dans le cas 
de Kummer ou d'Artin-Schreier (elle prend deux
fois moins de valeurs.)

\medskip 

{\bf Exemple }: on choisit $p=q=13$ et $d=7$ donc
le cofacteur est $m=2$.
On pose $D=2$ et on v{\'e}rifie que $D$ n'est pas un carr{\'e} dans
$\FF_{13}$.
On cherche $r=U+\sqrt 2 V$ tel que $U^2-2V^2=1$ et
$r$ soit d'ordre $p+1=14$ dans $\FF_{13}(\sqrt 2)^*$.
Par exemple $U=3$ et $V=2$ conviennent. La coordonn{\'e}e
$u$ de $3+2\sqrt 2$ est $u(r)=\frac{3}{2}=8$.
On est alors en mesure d'{\'e}crire le polyn{\^o}me 

$$A(X)=X^7+3X^5+10X^3+4X-8(7X^6+5X^4+6X^2+8).$$

En outre la formule \ref{eqn:frob} pr{\'e}dit l'action
du Frobenius.
On pose $t=[-m]r=[-2]r$ donc $u(t)=4$ et le Frobenius
op{\`e}re comme la translation par 
$t$ : 

$$X^p= \frac{4X+2}{X+4} \bmod A(X).$$

On a donc r{\'e}alis{\'e} un petit progr{\`e}s : d{\'e}sormais on sait traiter les
extensions de $\Fq$ dont le degr{\'e} $d$ divise $q+1$. Malheureusement
cette condition est aussi restrictive que celle impos{\'e}e par
la th{\'e}orie de Kummer.  Que faire si le degr{\'e} $d$ ne divise 
ni $q+1$ ni $q-1$ ?

Il faut diversifier les groupes alg{\'e}briques.
Les courbes elliptiques offrent une
alternative naturelle.

\section{Corps r{\'e}siduels sur les courbes elliptiques}\label{section:ell}

On revient {\`a} la d{\'e}marche de la section \ref{section:isog} en prenant
pour groupe alg{\'e}brique $\bG$ une courbe elliptique.

On consid{\`e}re un corps fini $\bK=\Fq$ dont on veut 
construire une extension de degr{\'e} $d$ avec $d$ premier {\`a} la
caract{\'e}ristique $p$ de $\Fq$. 

Soit donc $\bG=E$ une courbe elliptique ordinaire sur $\Fq$ 
et soit $\igot$ un id{\'e}al inversible 
de  l'anneau
d'endomorphismes $\End(E)$. On suppose
que $\igot$ divise  $\phi-1$ et que $\End(E)/\igot$ est cyclique
d'ordre $d$.
Donc  $E(\Fq)$ contient
un sous-groupe $T=\Ker \igot$ cyclique d'ordre $d$.

Soit $I : E \rightarrow F$ l'isog{\'e}nie cyclique de degr{\'e} $d$
et de noyau $T$. 
Le quotient $F(\Fq)/I(E(\Fq))$ est isomorphe {\`a} $T$.

Soit donc  $a$ dans $F(\Fq)$ tel que
$a\bmod I(E(\Fq))$  engendre ce quotient.

La fibre  $I^{-1}(a)$
est  un diviseur irr{\'e}ductible. Cela signifie que 
les $d$ points g{\'e}om{\'e}triques
au dessus de $a$ sont d{\'e}finis sur l'extension $\bL$ de degr{\'e}
$d$
 de $\bK$ et qu'ils sont conjugu{\'e}s entre eux par l'action
de Galois. On note
$B=I^{-1}(a)$ le diviseur premier correspondant.

Ainsi  $\bL$ est l'extension r{\'e}siduelle de $E$ en $B$.
Pour repr{\'e}senter un {\'e}l{\'e}ment de $\bL$ on se donne une fonction $f$ 
sur $E$ dont les p{\^o}les {\'e}vitent $B$ et on consid{\`e}re
l'{\'e}l{\'e}ment
$f\bmod B \in \bL$
appel{\'e}  r{\'e}sidu de $f$ en $B$. Soient $X$, $Y$, $Z$  des coordonn{\'e}es
projectives sur $E$.

Pour tout entier $k\ge 0$ on note $\cF_k$ l'ensemble des $\Fq$-fonctions
sur $E$ sans p{\^o}le en $B$ et de degr{\'e} $\le k$.

On note $P_k$ l'ensemble des {\'e}l{\'e}ments de $\bL$ correspondant

$$P_k=\{f\bmod B | f \in \cF_k     \}.$$

On a clairement (Riemann-Roch) 

$$\bK=P_0=P_1\subset P_2\subset \dots \subset P_{d}=\bL$$
\noindent et 

$$P_k\times P_l\subset P_{k+l}.$$

En outre il est clair que $\cF_k$ est invariant par $T$. Donc
$P_k$ est invariant par l'action de $\Gal(\bL/\bK)$.

Pour tester  si un {\'e}l{\'e}ment $z$ de $\bL$ est dans $P_k$ on cherche
une fonction $f$ dans $\cF_k$ telle que $f=z\pmod B$. C'est un probl{\`e}me
d'interpolation {\`a} peine plus difficile que dans les deux cas pr{\'e}c{\'e}dents
(polyn{\^o}mes pour Kummer et fractions rationnelles pour le tore).
Il suffit de chercher $f$ sous la forme $\frac{N}{D}$ o{\`u} $N$
et $D$ sont des formes homog{\`e}nes de  degr{\'e} $\lceil k/3\rceil +1$.
C'est encore un probl{\`e}me d'alg{\`e}bre lin{\'e}aire.

On peut prendre pour base de friabilit{\'e} l'ensemble des {\'e}l{\'e}ments $f\bmod
B$ de $P_\kappa$ avec $\kappa$ la borne de friabilit{\'e} choisie. 

Pour factoriser un {\'e}l{\'e}ment $z=f\bmod B$ de $\bL$ on d{\'e}compose  le
diviseur de $f$ en somme de diviseurs premiers et on esp{\`e}re que tous ces
diviseurs
ont un degr{\'e} $\le \kappa$.

On regarde maintenant quelles sont les conditions pour qu'existe
une courbe elliptique avec toutes les propri{\'e}t{\'e}s que nous 
avons requises.

 On veut
une courbe elliptique sur $\Fq$ de cardinal divisible par $d$. Donc
$q$ ne peut pas {\^e}tre trop petit. On doit avoir au moins

$$q+2\sqrt q+1>  d.$$

Pour simplifier on suppose que $d$ est impair
et admet un multiple  $D$ sans facteur
carr{\'e} tel que $D\not \equiv 1\bmod p$ et 

$$q+1-2\sqrt q < D < q+1+2\sqrt q.$$

Il existe alors 
une courbe elliptique ordinaire $E$  sur $\Fq$ de cardinal $D$ 
et de trace $q+1-D$. 

L'anneau $\ZZ[\phi]$  est int{\'e}gralement clos localement en chaque
premier impair divisant $D$, donc aussi  $\End(E)$.

L'id{\'e}al $(\phi -1)$ de $\End(E)$ admet un unique facteur $\igot$ de degr{\'e} $d$.
Le quotient $\End(E)/\igot$ est cyclique et $\igot$ est inversible
dans $\End(E)$.

La th{\'e}orie de la multiplication complexe (ou une simple
recherche exhaustive) permet de construire
la courbe $E$ une fois choisi $\phi$.

{\bf Exemple : } on choisit $p=q=11$ et 
$d=D=7$ donc $t=5$ et $\phi^2-5\phi+11=0$.
Le discriminant de $\ZZ[\phi]$ est
$-19$ donc $\End(E)=\ZZ[\phi]$. En particulier
 $\igot =(\phi-1)$ est inversible et son noyau $T$ est 
le groupe des points rationnels.

On consid{\`e}re l'isog{\'e}nie $I : E\rightarrow F$ de degr{\'e} $7$ obtenue en quotientant
$E$ par le groupe des points rationnels.

Pour tout  $a\in F(\FF_{11})$ non nul on sait que $B=I^{-1}(a)$ est  irr{\'e}ductible.

On trouve une {\'e}quation de $E$ : 

$$y^2 + xy = x^3 + 2x + 8.$$

\section{Les cribles en dimension deux}\label{section:JL}

Il existe une famille d'algorithmes pour la factorisation et le logarithme
discret, 
 appel{\'e}s crible alg{\'e}brique, crible du corps
des fonctions, etc., 
 qui reposent sur des calculs d'intersection sur une surface
({\'e}ventuellement
arithm{\'e}tique). Le principe de ces algorithmes est donn{\'e} en un seul dessin
sur la couverture de \cite{nfs}.

Nous illustrons ces id{\'e}es dans un
cadre un peu g{\'e}n{\'e}ral afin de pr{\'e}parer  l'exposition de notre 
construction dans la section \ref{section:carre}.

Au fil de notre exposition, nous illustrons  ce cadre
g{\'e}n{\'e}ral {\`a} travers l'un de ces algorithmes, d{\^u} {\`a} Joux et Lercier
\cite{JouxLercier}. 

Soit $\Fp$ le corps {\`a} $p$ {\'e}l{\'e}ments  avec $p$ premier.

On se donne une surface alg{\'e}brique projective lisse irr{\'e}ductible $\cS$ sur $\Fp$. 
Soient $\cA$ et $\cB$ deux courbes sur 
$\cS$.
Soit $\cI$ une sous-vari{\'e}t{\'e} irr{\'e}ductible de l'intersection $\cA\cap\cB$.
On suppose que $\cA$ et $\cB$ sont transverses en $\cI$ et 
on note $d$ le degr{\'e} de $\cI$.
Le corps r{\'e}siduel de $\cI$ est  donc $\Fp(\cI)=\Fq$ avec $q=p^d$.

Soit alors un pinceau (lin{\'e}aire ou du moins
alg{\'e}brique connexe) de  diviseurs $(D_\lambda)_{\lambda \in \Lambda}$ 
sur $\cS$. Ici $\Lambda$ est l'espace des param{\`e}tres.

On fixe un entier $\kappa$
et on recherche ({\`a} t{\^a}tons) des diviseurs $D_\lambda$ tels que
les deux diviseurs d'intersection $D\cap \cA$ et $D\cap \cB$
soient disjoints de $\cI$ et $\kappa$-friables (autrement dit, ils se
d{\'e}composent en somme de diviseurs de degr{\'e} $\le \kappa$.)

On consid{\`e}re la relation d'{\'e}quivalence $\equiv_\cI$ sur les diviseurs
de $\cS$, d{\'e}finie par $D\equiv_\cI 0$  si et seulement si $D$ est
le diviseur d'une fonction $f$ constante  modulo $\cI$.
Les classes d'{\'e}quivalence pour cette relation sont param{\'e}tr{\'e}es
par les points d'un groupe alg{\'e}brique $\Pic(\cS,\cI)$, extension
de $\Pic(\cS)$ par un tore $T_\cI$ de dimension $d-1$.

On d{\'e}finit de m{\^e}me les groupes alg{\'e}briques $\Pic(\cA,\cI)$ et
$\Pic(\cB,\cI)$ qui sont des jacobiennes g{\'e}n{\'e}ralis{\'e}es de $\cA$
et $\cB$ respectivement.

On a des morphismes $\Pic(\cS,\cI)\rightarrow \Pic(\cA,\cI)$ et
$\Pic(\cS,\cI)\rightarrow \Pic(\cB,\cI)$ qui induisent l'identit{\'e}
sur le tore $T_\cI$.

Soit $N$ un entier qui annule les trois groupes
$\Pic(\cS)(\Fp)$, $\Pic(\cA)(\Fp)$, et $\Pic(\cB)(\Fp)$.

Soient $\lambda$ et $\mu$ deux param{\`e}tres dans $\Lambda$
tels que  $D_\lambda\cap \cA$, $D_\mu \cap \cA$,
$D_\lambda\cap \cB$, et $D_\mu \cap \cB$ soient
friables. On suppose aussi que $D_\lambda$ et $D_\mu$ sont disjoints de $\cI$.

On {\'e}crit $D_\lambda\cap \cA = \sum A_i$,  $D_\mu\cap \cA = \sum B_j$, 
$D_\lambda\cap \cB = \sum C_k$,  $D_\mu\cap \cB = \sum D_l$ comme sommes de
diviseurs
sur $\cA$ ou $\cB$ de degr{\'e}s $\le \kappa $.

Le diviseur $D_\lambda-D_\mu$ est alg{\'e}briquement {\'e}quivalent {\`a} z{\'e}ro
et le diviseur $N(D_\lambda-D_\mu)$ est principal.

Soit $f$ une fonction sur $\cS$ de diviseur $N(D_\lambda-D_\mu)$.

On choisit un diviseur $X$ sur $\cA$ de degr{\'e} $1$ 
et un diviseur $Y$ sur $\cB$ de degr{\'e} $1$.

Pour tout $i$ soit  $\alpha_i$ une fonction sur $\cA$ 
de diviseur $N(A_i-\deg(A_i)X)$. 

Pour tout $j$ soit  $\beta_j$ une fonction sur $\cA$ 
de diviseur $N(B_j-\deg(B_j)X)$. 

Pour tout $k$ soit  $\gamma_k$ une fonction sur $\cB$ 
de diviseur $N(C_k-\deg(C_k)Y)$. 

Pour tout $l$ soit  $\delta_l$ une fonction sur $\cB$ 
de diviseur $N(D_l-\deg(D_l)Y)$.

On a

$$\frac{\prod_i\alpha_i}{\prod_j\beta_j}  = \frac{\prod_k\gamma_k
}{\prod_l\delta_l}\bmod \cI$$
\noindent ce qui produit une relation dans le groupe 
de $T_\cI(\Fp)=\Fq^*/\Fp^*$.

Par exemple Lercier et
Joux consid{\`e}rent $\cS=\PU \times \PU$. Pour {\'e}viter toute confusion on
note $\cC_1=\PU$ le premier facteur et $\cC_2=\PU$ le second facteur. 
Soit $O_1$ un point rationnel sur $\cC_1$ et soit $\cU_1=\cC_1-O_1$. Soit
$x$ une coordonn{\'e}e affine sur $\cU_1\sim \AU$. Soit de m{\^e}me $O_2$, $\cU_2$ et
$y$
une coordonn{\'e}e affine sur $\cU_2$.

Joux et Lercier choisissent pour $\cA$
l'adh{\'e}rence de Zariski dans $\cS$ de la courbe de $\cU_1\times \cU_2$
d'{\'e}quation $y=f(x)$ o{\`u} $f$ est un polyn{\^o}me de degr{\'e} $d_f$ dans $\Fp[x]$.
Pour $\cB$ ils choisissent l'adh{\'e}rence de Zariski dans $\cS$ de la courbe de $\cU_1\times \cU_2$
d'{\'e}quation $x=g(y)$ o{\`u} $g$ est un polyn{\^o}me de degr{\'e} $d_g$ dans $\Fp[y]$.

Le groupe de N{\'e}ron-Severi de $\cS$ est isomorphe {\`a} $\ZZ\times \ZZ$. La classe
d'{\'e}quivalence alg{\'e}brique d'un diviseur $D$ est donn{\'e}e par son bidegr{\'e}
$(\dx(D),\dy(D))$
avec $\dx(D)=D.(\cC_1\times O_2)$ et  $\dy(D)=D.( O_1\times
\cC_2)$. Et la forme d'intersection est donn{\'e}e par

$$D.E=\dx(E)\dy(D)+\dx(D)\dy(E).$$

Le bidegr{\'e} de $\cA$ est $(d_f,1)$ et celui de $\cB$ est $(1,d_g)$.
Ainsi $\cA.\cB=1+d_fd_g$ et l'intersection de $\cA$ et $\cB$ est form{\'e}e du point 
$O_1\times O_2$ et des $d_fd_g$ points de la forme $(\alpha,f(\alpha))$ o{\`u}
$\alpha$
est l'une des $d_fd_g$ racines de $g(f(x))-x$.

Soit alors $h(x)$ un facteur irr{\'e}ductible simple de ce dernier polyn{\^o}me  et soit $d$
son degr{\'e}.

On note $\cI$ la vari{\'e}t{\'e} de dimension $0$ et de degr{\'e} $d$ correspondante. Le
corps r{\'e}siduel $\Fp(\cI)$
est un  corps fini {\`a} $q$ {\'e}l{\'e}ments avec $q=p^d$.

Reste {\`a} construire  un pinceau de  diviseurs $(D_\lambda)_{\lambda \in \Lambda}$ 
sur $\cS$. Il est
naturel de consid{\'e}rer  l'ensemble $\Lambda$ des polyn{\^o}mes dans
$\Fp[x,y]$  de bidegr{\'e}  $(u_x,u_y)$ bien choisi.  Le diviseur $D_\lambda$
 correspondant
au polyn{\^o}me $\lambda$ 
 est la cl{\^o}ture de Zariski du lieu des z{\'e}ros de $\lambda$. Il 
a pour bidegr{\'e} $(u_x,u_y)$ lui aussi.

On fixe un entier $\kappa$
et on recherche ({\`a} t{\^a}tons) des diviseurs $D_\lambda$ tels que
les deux diviseurs d'intersection $D_\lambda\cap \cA$ et $D_\lambda \cap \cB$
soient disjoints de $\cI$ et $\kappa$-friables.

Par exemple, si $\lambda (x,y)$ est un polyn{\^o}me en $x$ et $y$, l'intersection
de $D_\lambda$ et de $\cA$ est de degr{\'e} $d_fu_y+u_x$. Sa partie affine est
 d{\'e}crite par les racines du polyn{\^o}me
$\lambda(x,f(x))=0$.

L'intersection
de $D_\lambda$ et de $\cB$ est de degr{\'e} $u_y+u_xd_g$. Sa partie affine est
 d{\'e}crite par les racines du polyn{\^o}me
$\lambda(g(y),y))=0$.

Il convient alors d'ajuster $u_x$ et $u_y$ en fonction de $p$ et $d$.

\section{Corps r{\'e}siduels sur des carr{\'e}s elliptiques}\label{section:carre}

Dans cette section on cherche {\`a} concilier la construction g{\'e}n{\'e}rique de la
section \ref{section:JL} et les id{\'e}es de la section 
\ref{section:ell}.

On demande que les  automorphismes de $\Fp(\cI)$ soient induits par des
automorphismes
de la surface $\cS$.

Soit donc
$E$ une courbe elliptique ordinaire sur $\Fp$ 
et soit $\igot$ un id{\'e}al inversible 
de  l'anneau
d'endomorphismes $\End(E)$. On suppose
que $\igot$ divise  $\phi-1$ et que $\End(E)/\igot$ est cyclique
d'ordre $d$.
Donc  $E(\Fq)$ contient
un sous-groupe $T=\Ker \igot$ cyclique d'ordre $d$.
Soit $I : E \rightarrow F$ l'isog{\'e}nie quotient par $\Ker \igot$ et
soit $J : F\rightarrow E$ telle que $\phi-1=J\circ I$.

On choisit pour surface $\cS$ le produit $E\times E$ et pour {\'e}viter toute
confusion
on note $E_1$ le premier facteur et $E_2$ le deuxi{\`e}me facteur. On note $O_1$
l'origine de $E_1$ et $O_2$ l'origine de $E_2$. 

Le groupe de N{\'e}ron-Severi est $\ZZ\times \ZZ \times \End(E)$. La classe $(d_1,d_2,\xi)$ d'un
diviseur
$D$ est form{\'e}e du bidegr{\'e}  et de l'isog{\'e}nie  induite par $D$. Plus pr{\'e}cis{\'e}ment
$d_1$
est le degr{\'e} d'intersection de $D$ et $E_1\times O_2$ et $d_2$ est le degr{\'e} de
$O_1\times E_2$ et $\xi : E_1\rightarrow E_2$.

Soient $\alpha$ et $\beta$ deux endomorphismes de $E$ et soient 
$a$ et $b$ deux points $\Fp$-rationnels sur $E$. 

Soit $\cA$ l'image r{\'e}ciproque de $a$ par l'application de $E\times E$
dans $E$ qui {\`a} $(P,Q)$ associe $\alpha(P)-Q$.

Soit $\cB$ l'image r{\'e}ciproque de $b$ par l'application de $E\times E$
dans $E$ qui {\`a} $(P,Q)$ associe $P-\beta(Q)$.

On suppose que $1-\beta\alpha=\phi-1$.  L'intersection 
 de $\cA$ et de $\cB$
est
form{\'e}e des couples $(P,Q)$ tels que $(\phi-1)(P)=b+\beta(a)$ et
$Q=\alpha(P)-a$.

On a choisi $a$ et $b$ de telle sorte
qu'il existe un point $c$ dans $F(\Fp)$ tel que $J(c)=b+\beta(a)$
et $c$ engendre $F(\Fp)/I(E(\Fp))$. L'intersection
de $\cA$ et $\cB$ contient alors une composante  $\cI$
irr{\'e}ductible 
de degr{\'e} $d$.

Soit maintenant $D$ un diviseur sur $\cS$ et $(d_1,d_2,\xi)$
sa classe dans le groupe de N{\'e}ron-Severi.

La classe de $\cA$ est $(\alpha\bar\alpha, 1,  \alpha)$ et
celle de $\cB$ est $(1,\beta\bar\beta, \bar \beta)$.

Le degr{\'e} d'intersection de $D$ et $\cA$ est donc

\begin{equation}\label{eq:int1}
D.\cA=d_x+d_y\alpha\bar\alpha  -\xi\bar \alpha-\bar \xi\alpha
\end{equation}

\noindent et de m{\^e}me
\begin{equation}\label{eq:int2}
D.\cB=d_x\beta\bar\beta+d_y  -\xi\bar \beta-\bar \xi\beta.
\end{equation}

On s'int{\'e}resse particuli{\`e}rement au cas o{\`u}
 les normes de $\alpha$ et $\beta$
sont de tailles comparables (soit la racine carr{\'e}e de la norme 
de $\phi -2$).

On obtient alors des performances comparables {\`a} celles de la section
\ref{section:JL} mais avec un avantage : les bases de friabilit{\'e}s sur
$\cA$ et sur $\cB$ sont invariantes par action de Galois.

Soit en effet $f$ une  fonction de degr{\'e} $\le \kappa$ sur
$\cA$. Un point de $\cA$ est un couple $(P,Q)$ tel que
$Q=\alpha(P)-a$. On l'identifie donc {\`a} sa coordonn{\'e}e $P$
et on voit $f$ comme une fonction sur  $E_1$. Supposons en outre
que $(P,Q)$ est dans $\cI$. Alors $f(P,Q)=f(P)$ est 
un {\'e}l{\'e}ment de la base de friabilit{\'e} sur $\cA$. On observe
alors que $f(P)^p=f(\phi(P))=f(P+t)$ o{\`u} $t$ est un {\'e}l{\'e}ment du
noyau $T$ de $\igot$. Donc $f(P)^p$ est la valeur en $P$ de
$f\circ\tau_{t}$ avec $\tau_t : E_1\rightarrow E_1$
la translation par $t$.
Comme $f\circ \tau_t$ est une fonction de m{\^e}me degr{\'e} que $f$,
sa valeur en $P$ est encore un {\'e}l{\'e}ment de la base de friabilit{\'e}.

On peut ainsi diviser par $d$ la taille de la base de friabilit{\'e} 
sur $\cA$ et aussi sur $\cB$.

On choisit de petites valeurs de $(d_x,d_y,\xi)$ en pr{\'e}f{\'e}rant celles
qui minimisent les expressions \ref{eq:int1} et \ref{eq:int2}.
On demande que  $d_x\ge 1$, $d_y\ge 1$ et

\begin{equation}\label{eq:ine}
d_xd_y\ge \xi\bar\xi+1. 
\end{equation}

On a  donc une
classe d'{\'e}quivalence alg{\'e}brique 
$\cgot = (d_x,d_y,\xi)$.

Et on cherche  les diviseurs
effectifs de cette classe.

On a  un point  $O_1$ sur  $E_1$  et  un point $O_2$ sur $E_2$.

Le graphe $\cG = \{(P,Q) | Q=-\xi (P)  \}$ de $-\xi : E_1\rightarrow E_2$ est un
diviseur de la classe $(\xi\bar\xi , 1, -\xi)$ donc 
$\cH=-\cG+(d_x+\xi\bar\xi )O_1\times E_2+(d_y+1)E_1\times O_2$ 
est dans $\cgot$.

On  calcule l'espace lin{\'e}aire 
$\cL(-\cG+(d_x+\xi\bar\xi )O_1\times
E_2+(d_y+1)E_1\times O_2)$ en utilisant la suite de restriction

$$0\rightarrow \cL_S(-\cG+(d_x+\xi\bar\xi)O_1\times E_2
+(d_y+1)E_1\times O_2)
\rightarrow \cL_{E_1}((d_x+\xi\bar\xi )O_1)\otimes \cL_{E_2}((d_y+1) O_2)
\rightarrow \cL_\cG(\Delta)$$
\noindent o{\`u} $\Delta$ est le diviseur sur $\cG$ donn{\'e}
par l'intersection avec 

$$(d_x+\xi\bar\xi)O_1\times E_2+(d_y+1)E_1\times O_2.$$

Ce diviseur est de degr{\'e} $d_x+\xi\bar\xi +(d_y+1)\xi\bar\xi$ donc la
dimension
du terme de droite dans la suite ci-dessus est {\'e}gale {\`a} ce nombre.

D'autre part, le terme du milieu est de dimension
 $(d_x+\xi\bar\xi)(d_y+1)$
qui est strictement sup{\'e}rieur {\`a} la dimension du terme de droite ({\`a} cause
de l'in{\'e}galit{\'e} \ref{eq:ine}).

Donc l'espace lin{\'e}aire de gauche est non nul et la classe est
effective.

Ainsi, la condition num{\'e}rique \ref{eq:ine} est un crit{\`e}re suffisant
d'effectivit{\'e}.

En pratique on calcule
une base de 
$\cL_{E_1}((d_x+\xi\bar\xi)O_1)$
et une base de $\cL_{E_2}((d_y+1) O_2)$ et on multiplie
les deux bases (on prend tous les produits form{\'e}s d'un {\'e}l{\'e}ment
de la premi{\`e}re base et d'un {\'e}l{\'e}ment de la deuxi{\`e}me).

On s{\'e}lectionne un nombre suffisant  (plus de 
$d_x+\xi\bar\xi+(d_y+1)\xi\bar\xi$) de points 
$(A_i)_i$ sur
$\cG$ et on {\'e}value toutes les fonctions en ces points. 
Un calcul d'alg{\`e}bre  lin{\'e}aire donne une base
de l'ensemble des fonctions qui s'annulent
en tous ces points, donc aussi le long de $\cG$.

Pour chaque fonction $\phi$ dans cet espace, le diviseur
des z{\'e}ros de $\phi$ contient $\cG$ et la diff{\'e}rence $(\phi)_0
-\cG$ est un  diviseur effectif dans la classe d'{\'e}quivalence lin{\'e}aire de $\cH$.

On a donc construit une  classe d'{\'e}quivalence lin{\'e}aire dans $\cgot$.
Pour construire les autres classes d'{\'e}quivalences lin{\'e}aires de $\cgot$ on
note que  $E\times E$ est sa propre  vari{\'e}t{\'e} de Picard. Il suffit
donc de remplacer $\cH$ dans  le calcul pr{\'e}c{\'e}dent par
$\cH+E_1\times Z_2-E_1\times O_2+Z_1\times E_2-O_1\times E_2$ 
o{\`u} $Z_1$ et $Z_2$ parcourent $E_1(\Fp)$ et $E_2(\Fp)$ respectivement.

\section{G{\'e}n{\'e}ralisation et limites ?}\label{section:conclusion}

La construction de la section \ref{section:carre} peut et doit  {\^e}tre
g{\'e}n{\'e}ralis{\'e}e.

Soit encore 
$E$ une courbe elliptique ordinaire sur $\Fp$ 
et soit $\igot$ un id{\'e}al inversible 
de  l'anneau
d'endomorphismes $\End(E)$. On suppose
que $\igot$ divise  $\phi-1$ et que $\End(E)/\igot$ est cyclique
d'ordre $d$. On note $F$ le quotient de $E$ par le noyau de
$\igot$ et $I : E \rightarrow F$ l'isog{\'e}nie quotient.

L'entier $d$ appartient {\`a} l'id{\'e}al $\igot$. Soient $u$ et veux deux {\'e}l{\'e}ments
de $\igot$ tels que $d=u+v$ et $(u)=\igot\agot_1\bgot_1$
et $(v)=\igot\agot_2\bgot_2$ o{\`u} $\agot_1$, $\bgot_1$, $\agot_2$, $\bgot_2$
sont des id{\'e}aux inversibles de $\End(E)$.

On en d{\'e}duit l'existence de deux courbes $E_1$ et $E_2$ et de quatre
isog{\'e}nies $\alpha_1$, $\beta_1$, $\alpha_2$, $\beta_2$, telles que
$\beta_1\alpha_1+\beta_2\alpha_2=I$. On repr{\'e}sente ci-dessous
ces trois isog{\'e}nies de $E$ vers $F$ 

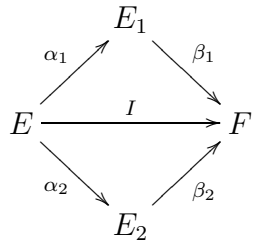
\begin{eqnarray*}\xymatrix{
 & E_1 \ar@{->}^{\beta_1}[dr]   &   \\
E \ar@{->}^{I}[rr]\ar@{->}^{\alpha_1}[ur]\ar@{->}_{\alpha_2}[dr]&  &F\\
& E_2  \ar@{->}_{\beta_2}[ur]&
}
\end{eqnarray*}
\noindent 

On choisit $\cS=E_1\times E_2$. Pour $\cA$ on choisit l'image
de  $(\alpha_1,\alpha_2) : E\rightarrow \cS$. Pour
$\cB$ on choisit l'image inverse de $f$ par $\beta_1+
\beta_2 : \cS\rightarrow F$ o{\`u} $f$ est un g{\'e}n{\'e}rateur de  $F(\Fp)/I(E(\Fp))$.

L'intersection de $\cA$  et $\cB$ est donc l'image par $\alpha_1\times
\alpha_2$  de $I^{-1}(f)\subset E$.

On choisit $u$ et $v$ de telle sorte que $\agot_1$, $\bgot_1$,
$\agot_2$,  et $\bgot_2$, aient des normes proches de la racine
carr{\'e}e de $d$.

Cette construction est utile lorsque la norme de $\igot$ est 
beaucoup plus petite
que celle de $\phi-1$.

Nous avons donc r{\'e}ussi {\`a} construire des bases invariantes de friabilit{\'e}
pour un grand nombre de corps finis. Nos constructions vont au del{\`a} des
 th{\'e}ories de Kummer et Artin-Schreier. Elles sont efficaces si le degr{\'e}
$d$ du corps est inf{\'e}rieur {\`a}  $4\sqrt q$ ou contenu dans l'intervalle
$]q+1-2\sqrt q,q+1+2\sqrt q [$.

\bibliography{smooth1}
\end{document}